\documentstyle[amsfonts,12pt]{article}
\newcommand{\hao}{\ha_1}
\newcommand{\goo}{\go_1}
\newcommand{\restriction}{{|}} 

\def\rest{\mathord{\restriction}}

\newcommand{\open}{\Bbb}
 
\def\deq{\mathop=\limits^{\rm def}} 

\newlength{\labparwidth}
\newcommand{\labpar}[2]{$$\parbox{\labparwidth}{#2}\leqno(#1)$$}

\newcommand{\ann}{\mbox{\rm Ann}}

\newcommand{\supp}{\mbox{\rm supp}}
\newcommand{\se}{\subseteq}
\newcommand{\set}[2]{\{#1 \colon #2\}} 
\newcommand{\qed}{$\Box$\par\medskip}
\newcommand{\dmd}{\diamondsuit}

\newcommand{\ga}{\alpha}
\newcommand{\gb}{\beta}
\newcommand{\grg}{\gamma}
\newcommand{\gd}{\delta}

\newcommand{\gz}{\zeta}

\newcommand{\gs}{\sigma}
\newcommand{\gt}{\tau}

\newcommand{\go}{\omega}
 
\newcommand{\ha}{\aleph}
 
\def\psupp{\mathop{\hbox{p-supp}}}
\begin{document}
\thispagestyle{empty}
\vspace*{0.2in}
\begin{center}
{\large \bf  ON A CONJECTURE REGARDING\\
 NON-STANDARD UNISERIAL MODULES}
\end{center}

\vspace{0.75in}

\begin{center}
 Paul C. Eklof \footnote{Thanks to Rutgers University for its support
of this research through its funding of the first author's visits to
Rutgers.} \\University of
California, Irvine\\  and\\
Saharon Shelah \footnote{Partially supported by Basic Research Fund,
Israeli Academy of Sciences. Publication 422.} \\ Hebrew University
\\ and Rutgers University
\end{center}
\vspace{1.0in}

\noindent
{\bf Abstract:} {\small We consider the question of which valuation domains (of 
cardinality $\aleph _1)$ have non-standard uniserial modules. We show that a 
criterion conjectured by Osofsky is independent of ZFC + GCH.}

\vspace{1.25 in}

\noindent
{\small {\it 1991 Mathematics Subject Classification}. Primary 13L05,
03E35,
13C05; Secondary 03E75, 13A18.

}

\newpage

\noindent
{\bf Introduction.} 

\medskip

The story of the study of non-standard uniserial modules is a long and 
interesting one, which we outline only briefly here; the reader may consult 
[BS1], [FS] or [O2] for more information. The existence of
non-standard uniserial modules over some valuation domains was first proved
 by the second author [Sh], by forcing a model of ZFC
with a non-standard uniserial module, and then 
using a completeness theorem for stationary logic to show that such examples 
exist in all models of ZFC, i.e., the existence is a consequence of ZFC. Then
Fuchs and Salce [FS] constructed non-standard uniserial modules using the 
diamond principle, which is consistent with ZFC.
Fuchs noticed that a non-standard divisible uniserial module could be
used to construct an affirmative  answer to Kaplansky's question whether there is a 
valuation ring which is not the quotient of a valuation domain.
 Fuchs and Shelah [FSh] used 
the compactness theorem for the logic with quantification over branches to give
another proof that divisible non-standard uniserials exist in all models of 
ZFC. The first author [E] elaborated on this argument to give a general 
transfer principle that, for example, showed that various different classes of 
non-standard uniserials constructed by Bazzoni and Salce [BS1, 2, 3] using 
$\diamondsuit $ existed in all models of ZFC.  

Now this use of theorems from mathematical logic did not bother the second 
author; he was happy if algebraists had to learn logic! However, in this view 
he was probably a singleton among those interested in the problem. Aside from 
eliminating the methods of logic from the proof, the hope was that a more 
explicit proof would give more information, e.g. an understanding of exactly 
which valuation domains had non-standard uniserial modules. Osofsky [O1, 2] was
the first to give a concrete construction (in ZFC) of valuation domains with 
non-standard uniserials. She ventured a conjecture [O1] as to which valuation 
domains had non-standard uniserials; we will state it below in a modified form.
It is a principal result of this paper that, while this conjecture is correct 
assuming  V = L (the Axiom of Constructibility), it fails in some other models of ZFC + GCH. Thus, it seems,
logic cannot be expelled from the subject!  

It is well-known that if $R$ is ``complete" in some sense (e.g. almost 
maximal), then there are no non-standard uniserial $R$-modules. Here we 
introduce a couple of invariants (``Gamma invariants") of a ring
whose values are stationary subsets of $\omega _1$ --- or, more
precisely, equivalence classes of stationary subsets under the
equivalence relation of equality on a closed unbounded subset.
These invariants can be viewed as measuring (the lack of) ``completeness"
of the ring. Some
results about implications between the values of these invariants and the 
existence of non-standard uniserials are proved in ZFC, but other implications 
are shown to be independent of ZFC. 

We would like to thank Luigi Salce and Silvana Bazzoni for helpful comments.

\bigskip

\bigskip

\noindent
{\bf Preliminaries.} 

\medskip
 An $R$-module $U$ is called {\em uniserial} if its
submodules are totally ordered by inclusion. An integral domain $R$
 is a valuation domain if it is uniserial as an $R$-module.
Throughout, $R$ will denote a valuation domain and $Q$  the quotient field of
$R$; we assume $ Q \neq R$. The uniserial modules we consider will always be generated by (at most)
$\aleph _1$ elements, and most of the time $R$
will be of cardinality $\aleph _1$. $R^{*}$ denotes the set of units of $R$. The 
residue field of $R$ is defined to be $R/P$, where $P$ is the maximal ideal of 
$R$.

If $J$ and $A$ are $R$-submodules of $Q$ with $A \subseteq J$, then
$J/A$ is a uniserial $R$-module, which is said to be  {\it standard}.
 A uniserial
$R$-module $U$ is said to be {\it non-standard} if it is not isomorphic
to a standard uniserial.

Given a uniserial module $U$, and a non-zero
element, $a$, of $U$, let $\ann(a) = \set{r \in R}{ra = 0}$ and
let $D(a) = \cup \{r^{-1}R \colon r$ divides a
in $U\}$. We say $U$ is  of {\em type} $J/A$ if $J/A \cong D(a)/\ann(a)$.
This is well-defined in that if $b$ is another non-zero element of $U$,
then $D(a)/\ann(a) \cong D(b)/\ann(b)$. For example,
$U$ has type $Q/R$
if and only if $U$ is divisible torsion and  the annihilator
ideal of every non-zero element of
$U$ is principal. (But notice that there is no $a \in U$ with $\ann(a) =
R$.) It is not hard to
see that if $U$ has type $J/A$, then
 $U$ is standard if and only if it is isomorphic to
$J/A$.

Suppose $U$ has type $J/A$ and there exists 
 $a \neq 0$ in $U$ such that $A = \ann(a)$ and $J = D(A)$.
 Then
$$
J = \cup _{\nu <\omega _1}r^{-1}_\nu R
\leqno(*)
$$
for some sequence of elements $\{r_\nu \colon \nu  < \omega _1\}$ such that for all 
$\mu  < \nu $, $r_\mu |r_\nu $.  If $J$ is countably generated, then $U$ is 
standard, so generally we will be assuming that $J$ is not countably generated;
then it has a set of generators as in ($\ast$), where, furthermore, $r_\nu $ does 
not divide $r_\mu $ if $\mu  < \nu $, i.e. $r^{-1}_\nu r_\mu $ is not a member
of $R$. Note that if $a_\nu  \in  U$ such that $r_\nu a_\nu  = a$,
then $\ann(a_\nu) = r_\nu A$, so $a_\nu R \
(\subseteq  U)$ is isomorphic to $R/r_\nu A = R/(A : r^{-1}_\nu R)$. If
$\delta  \in  \lim (\omega _1)$ and 
$$
J_\delta 
\deq \cup _{\nu <\delta } r^{-1}_\nu R
\leqno(\ast \ast)
$$
then the submodule of $U$ generated by $\{a_\nu \colon \nu  < \delta \}$ is 
a module over $R/(A : J_\delta ) = R/\cap _{\nu <\delta }r_\nu A.$

\bigskip

\noindent
{\bf Definition 1.} A subset $C$ of $\go_1$ is called a {\em cub} ---
short for closed unbounded set --- if $\sup C = \go_1$ and for all
 $Y \subseteq  C$,
$\sup Y \in  \go_1 $ implies  $\sup Y \in  C$.
 Call two subsets, $S_1$ and $S_2$, of $\go_1$
equivalent iff there is a cub $C$ such that $S_1 \cap C = S_2 \cap C$.
Let $\tilde{S}$ denote the equivalence class of $S$. The inclusion
relation induces a partial
order on the set, $D(\goo)$, of equivalence classes, i.e.,
$\tilde{S_1} \leq \tilde{S_2}$ if and only if there is a cub $C$ such
that $S_1 \cap C \se S_2 \cap C$. In fact, this induces a Boolean
algebra structure on $D(\goo)$, with least element, $0$, the equivalence
class of sets disjoint from a cub; and greatest element, $1$, the
equivalence class of sets containing a cub. We say $S$ is {\em
stationary} if $\tilde{S} \neq 0$, i.e., for every cub $C$, $C \cap S
\neq \emptyset$. We say $S$ is {\em co-stationary} if $\goo \setminus S$
is stationary.

Given $R$ and a type $J/A$, define $\Gamma _R(J/A)$ to be
$\tilde{S}$, where
$$
S = \{\delta  \in  \lim (\omega _1)\colon R/(A : J_\delta )\hbox{ is not complete}\}
$$

\noindent
where the topology on $R/(A : J_\delta )$ is the metrizable linear topology 
with a basis of neighborhoods of 0 given by the submodules $r_\nu A\  (\nu  < 
\delta ).$ 
$\Gamma _R(J/A)$ is well-defined in the sense that it does not depend on the
choice of the sequence $\{r_\nu \colon \nu  < \omega _1\}$, because of the 
following, which is proved by a standard argument (cf. [EM; pp 85f]):

\bigskip

\noindent
{\bf LEMMA 2.}   {\it If $J = \cup _{\nu <\omega _1}r^{-1}_\nu R$ and also $J =
\cup _{\nu <\omega _1}s^{-1}_\nu R$, then $\{\delta \colon \cap _{\sigma <\delta } 
r_\sigma A = \cap _{\sigma <\delta } s_\sigma A\}$ is a cub.} \qed

\bigskip

Note that $\Gamma _R(J/A) =  0$ if $J$ is countably generated.

If $J/A \cong J'/A'$,
 then $\Gamma _R(J/A) = \Gamma _R(J'/A')$, but in general
$\Gamma _R(J/A)$ and $\Gamma _R(J'/A')$ will be different for different types 
$J/A$ and $J'/A'$. For example, if $R$, $A$ and $J$ are as in Example 2 of [BS1; 
p. 302], then $\Gamma _R(J/A) = 1$, but $\Gamma _R(Q/A) = 0$ (since $Q$ is 
countably generated). Note that if $R$ is almost maximal, then
$\Gamma_R(J/A) = 0$ for all types $J/A$.

\bigskip

In proving results about the existence or non-existence of non-standard 
uniserial modules it will be convenient to make use of the formulation of the problem given in
[BS1], especially Lemma 1.2 and Test Lemma 1.3. Thus any uniserial module $U$
of type $J/A$ is described up to isomorphism by a family of units,
$\{e^\tau _\sigma \colon \sigma  < \tau  < \omega _1\}$ such that 
$$
e^\delta _\tau e^\tau _\sigma  - e^\delta _\sigma  \in  r_\sigma A
\leqno(\dag)
$$
for all $\sigma  < \tau  < \delta  < \omega _1.$ 

Moreover, if $U$ is given by (\dag ), then $U$ is standard if and only if there
exists a family $\{c_\sigma \colon \sigma  < \omega _1\}$ of units of $R$ such that
$$
c_\tau  - e^\tau _\sigma c_\sigma  \in  r_\sigma A
\leqno{(\dag\dag)}
$$
for all $\sigma  < \tau  < \omega _1.$

Salce has pointed out that by results of [BFS], the question of the
existence of a non-standard uniserial $R$-module of type $J/A$ 
can be reduced to the question of  the existence of a non-standard uniserial of type $K/R$
for an appropriate $K$.

\bigskip

\bigskip

\noindent
{\bf The conjecture.} 

\medskip

We can paraphrase the conjecture of Osofsky in [O1], for valuation domains $R$ 
of cardinality $\aleph _1$, as follows:  
\begin{quotation}
\noindent
$R$\ {\it has a non-standard uniserial of type $J/A$ if and only if 
$\Gamma _R(J/A) = 1.$} 
\end{quotation}
\noindent
Now, assuming CH (and the Axiom of Choice), it is possible to construct a valuation domain $R$ of
cardinality $\hao$ such that $\Gamma_{R}(J/A) \neq 0, 1$ (see Definition
10 and the proof of Theorem 11).
Thus the condition $``\Gamma _R(J/A) = 1"$ cannot be necessary,
because of the following result which implies that it is consistent with ZFC 
that there is a non-standard uniserial $R$-module of type $J/A$
for such an $R$. Note that $\diamondsuit _{\omega _1}(S)$, for any
stationary $S$, is consistent with ZFC; in fact it is implied by  V = L
(cf. [D; p. 139]).

\bigskip

\noindent
{\bf PROPOSITION 3.}   {\it Assume $R$ has cardinality $\aleph _1$}.

(i) {\it If
$\Gamma _R(J/A) = \tilde{S}$, then $\diamondsuit _{\omega _1}(S)$ implies that there
exists a non-standard uniserial $R$-module of type $J/A$.}

(ii) {\it If $\Gamma(J/A) = 1$, then CH implies that there
exists a non-standard uniserial $R$-module of type $J/A$.}

\bigskip

 Proposition 3(i) is essentially due to Fuchs and Salce [FS]; it can
be proved, for example, by a straightforward generalization of the proof of
Theorem 1.4 of [BS1], but for later purposes we will sketch the proof of
a somewhat
stronger result here. We begin with a lemma which, in one or form or
another, is a staple of the subject; its statement and proof, in
our chosen notation, we include here for the sake of completeness.

\bigskip

\noindent
{\bf LEMMA 4.}   {\it Suppose $\{e^\tau _\sigma \colon \sigma  < \tau  < \delta \}$ 
is a family of units satisfying (}\dag {\it ) and $R/\cap _{\nu <\delta } 
r_\nu A$ is not complete. Then there are units $e^\delta _{\sigma ,j}
\ (\sigma  
< \delta $; $j = 0,1)$ such that}

$$
e^\delta _{\tau,j} e^\tau _\sigma  - e^\delta _{\sigma,j}  \in  r_\sigma A
$$

\noindent
{\it for all $\gs < \gt$ and $j \in \{0, 1\}$, but there is no sequence
$\{c_\sigma \colon \sigma  < \delta \}$ for which there are $c_{\delta, j} \in  R^{*}\ (j= 0,1)$ 
such that (}\dag \dag {\it ) holds for all $\sigma  < \tau  < \delta $ and 
moreover}
$$
c_{\delta, j} - e^\delta _{\sigma ,j}c_\sigma  \in  r_\sigma A
$$

\noindent
{\it for all $\sigma  < \delta $ and $j \in \{0, 1\}$}.

\medskip

\noindent
{\sc Proof}.  Let $e^\delta _{\sigma ,0}$ be any family of units such that 
$e^\delta _{\sigma ,0} - e^\delta _{\tau ,0}e^\tau _\sigma  \in  r_\sigma A$ 
for all $\sigma  < \delta $. Let $\langle v_\sigma \colon \sigma  < \delta \rangle $
represent a Cauchy sequence in $R/\cap _{\nu <\delta } r_\nu A$ which does not 
have a limit in $R/\cap _{\nu <\delta } r_\nu A$. Let $e^\delta _{\sigma ,1} = v_\sigma e^\delta _{\sigma ,0}$. 
Then it is easy to check that $e^\delta _{\sigma ,1} - 
e^\delta _{\tau ,1}e^\tau _\sigma  \in  r_\sigma A$ for all $\sigma  < 
\delta $. Moreover, if the conclusion of the lemma is contradicted by 
$\{c_\sigma \colon \sigma  < \delta \}$ and $c_{\gd,0}$, $c_{\gd,1}$,
then $c_{\gd,1}c^{-1}_{\gd,0}$ is a
limit of $\langle v_\sigma \colon \sigma  < \delta \rangle $, which is impossible.
\qed

\bigskip

\noindent
{\sc Proof} of Proposition 3\ (sketch). In fact we will prove that,
under the hypotheses on $R$, $\Phi _{\omega _1}(S)$ implies that there exists a
non-standard uniserial $R$-module of type $J/A$.
(This generalizes the result in [FrG].)  Here
 $\Phi _{\omega _1}(S)$ is the weak diamond principle which is implied
by $\dmd_{\go_1}(S)$ (see, for example, [EM; VI.1.6]). In [DSh] it is proved that CH implies
$\Phi _{\omega _1}(S)$, so part (ii) follows.

 Write $R$ as the union, $R = \cup _{\nu <\omega _1}R_\nu $, of
a continuous chain of countable subsets. For each $\delta  \in  S$, define a 
function $P_\delta $ on pairs $(e$, $c)$ where $e = \{e^\tau _\sigma \colon \sigma  
< \tau  < \delta \} \subseteq  R_\delta $, $c = \{c_\sigma \colon \sigma  < 
\delta \} \subseteq  R_\delta $ as follows. If $e$ satisfies (\dag ), let 
$\{e^\delta _{\sigma ,j}\colon \sigma  < \delta \} \  (j = 0,1)$ be as in Lemma 4; if $e$ and $c$ satisfy (\dag \dag ) and there is $c_0 \in  R^{*}$ such that $c_0 -
e^\delta _{\sigma ,0}c_\sigma  \in  r_\sigma A$  for all $\sigma  < \delta $, 
let $P_\delta (e$, $c) = 1$; in all other cases let $P_\delta (e$, $c) = 0$. 
Let $\rho \colon S \rightarrow  2$ be a weak diamond function (which predicts the 
values of the $P_\delta )$. Now define the $e^\tau _\sigma $ by induction. The 
crucial case is when $\delta  \in  S$ and $e = \{e^\tau _\sigma \colon \sigma  < 
\tau  < \delta \}$ is contained in $R_\delta$; in this case let 
$e^\delta _\sigma  = e^\delta _{\sigma ,\rho (\delta )}$.  
After completing the inductive construction of the $e^\tau_\sigma$
for all $ \sigma < \tau < \omega_1$, notice that for any sequence $\{c_\sigma
 \colon \sigma < \omega_1 \}$, there is a cub $C$ such that for all
$\delta \in C$, $e^\tau_\sigma \in R_\delta$ and $c_\sigma
\in R_\delta$ for all $\sigma < \tau < \delta$.
Then the construction and the properties of the weak diamond function
imply that there is no sequence $\{ c_\sigma \colon \sigma < \omega_1
\} $ satisfying (\dag\dag).
\qed

\bigskip

The condition $\Gamma _R(J/A) \neq  0$ is certainly
necessary for the existence of a non-standard
uniserial $R$-module of type $J/A$, as the following easy lemma shows.

\bigskip

\noindent
{\bf LEMMA 5.}   {\it If $\Gamma _R(J/A) = 0$, then every uniserial $R$-module 
of type $J/A$ is standard}.

\medskip

\noindent
{\sc Proof}.  Without loss of generality, we can assume that for all limit $\delta $,
$R/(A : J_\delta )$ is complete. Given $\{e^\tau _\sigma \colon \sigma  < \tau  < 
\omega _1\}$ as in (\dag ), we define the units $c_\sigma $ satisfying 
(\dag \dag ) by induction on $\sigma $. If $c_\sigma $ has been defined, let 
$c_{\sigma +1} = e^{\sigma +1}_\sigma c_\sigma $. For a limit ordinal 
$\delta $, consider the Cauchy sequence in $R/(A : J_\delta )$ represented by 
the elements $\langle e^\delta _\sigma c_\sigma \colon \sigma  < \delta \rangle $; 
this sequence has a limit which is represented by a unit of $R$, which we call 
$c_\delta $. \qed

\bigskip

Thus we end up with the following reformulated version of the conjecture:

\labpar{\$}{{\it For any valuation domain $R$ of cardinality $\aleph _1$, $R$ 
has a non-standard uniserial module of type $J/A$ if and only if 
$\Gamma _R(J/A) \neq  0.$}}

\bigskip

\noindent
{\bf THEOREM 6.}   {\it Assuming  V = L, this conjecture is true.}

\medskip

\noindent
{\sc Proof}.  In fact, this follows from Proposition 3 and Lemma 5 since G\"odel's 
Axiom of Constructibility,  V = L, implies $\diamondsuit _{\omega _1}(S)$ for 
all stationary $S$.  \qed

\bigskip

Thus the conjecture is consistent with ZFC, and in fact with ZFC +
GCH, since V = L implies GCH.

\bigskip

\bigskip

\noindent
{\bf Another Gamma Invariant.} 

\medskip

Let\ $R$ be a valuation domain of cardinality $\aleph _1$ and $J$, $A$ as before.
For any limit ordinal $\delta  < \omega _1$, let 

$${\cal T}^\delta _{J/A} = \{\langle u_\sigma \colon \sigma  < \delta \rangle \colon 
\forall \sigma  < \tau  < \delta (u_\sigma  \in  R^{*},\  \hbox{and } u_\tau  - u_\sigma 
\in  r_\sigma A)\};$$

\noindent
that is, ${\cal T}^\delta _{J/A}$ consists of sequences of units which are 
Cauchy in the metrizable topology on $R/(A: J_\delta )$. Let ${\cal L}^\delta _{J/A}$ consist of those 
members of ${\cal T}^\delta _{J/A}$ which have limits in $R$, i.e. 

$${\cal L}^\delta _{J/A}= \{\langle u_\sigma \colon \sigma  < \delta \rangle  \in  
{\cal T}^\delta _{J/A} \colon \exists u_\delta  \in  R^{*} \hbox{
s.t. }\forall \sigma  < 
\delta (u_\delta  - u_\sigma  \in  r_\sigma A)\}.$$ 
Note that $\Gamma _R(J/A) = \tilde{S}$ where
$$
S = \{\delta  \in  \lim (\omega _1)\colon {\cal T}^\delta _{J/A} \neq
{\cal L}^\delta _{J/A}\}.
$$

We will be making use of an $\omega _1${\it -filtration of $R$ by subrings}, by
which we mean an increasing chain $\{N_\alpha \colon \alpha  \in  \omega _1\}$ of 
countable subrings of $R$ such that $R = \cup _{\alpha \in \omega _1} 
N_\alpha $; for limit $\alpha $, $N_\alpha  = \cup _{\beta <\alpha } N_\beta $;
and for all $\alpha $, $R^{*} \cap  N_{\alpha} = N^*_{\alpha}$, and $r_\nu  \in  N_\alpha  
\Longleftrightarrow  \nu  < \alpha $ [where the $r^{-1}_\nu $ generate $J$ 
as in ($\ast$)].

\medskip

\noindent
{\bf Definition.}  $\Gamma '_R(J/A) =
\tilde{E}'$ where  
\begin{quotation}
\noindent
$E' =  \{\delta  \in  \lim (\omega _1)\colon \exists \langle u_\sigma \colon \sigma  < 
\delta \rangle  \in  {\cal T}^\delta _{J/A} \hbox{ s.t. } \forall \langle f_\sigma \colon 
\sigma  < \delta \rangle  \in  {\cal L}^\delta _{J/A} \exists  \sigma < 
\delta \ \hbox{s.t. }u_\sigma f_\sigma  \notin  N_\delta \}$\\ 
$ =  
\{\delta  \in  \lim (\omega _1)\colon \exists \langle u_\sigma \colon \sigma  < 
\delta \rangle  \in  {\cal T}^\delta _{J/A} \hbox{ s.t. } \forall f \in  R^{*} \exists  
\sigma  < \delta \ \hbox{ s.t. }   
u_\sigma f \notin  N_\delta \hbox{ mod } r_\sigma A\}.$
\end{quotation}

As before, by standard methods it can be shown that the definition does not 
depend on the choice of $\{r_\nu \colon \nu  < \omega _1\}$ or of $\{N_\alpha \colon 
\alpha  < \omega _1\}$. (Here it is important that the $N_\ga$ are
countable.) Notice that $\Gamma '_R(J/A) \leq  \Gamma _R(J/A)$
since if ${\cal T}^\delta _{J/A} = {\cal L}^\delta _{J/A}$, then we can let
$\langle f_\sigma \colon \sigma  < \delta \rangle $  be $\langle u^{-1}_\sigma \colon 
\sigma  < \delta \rangle $ to show $\delta \notin E'$. We will prove below that the condition 
$\Gamma _R(J/A) \neq  0$ is not (provably from ZFC) sufficient for the 
existence of a non-standard uniserial of type $J/A$; however, we do have the 
following theorem of ZFC.

\bigskip                   

\noindent
{\bf THEOREM 7.}   {\it If $\Gamma '_R(J/A) \neq  0$, then there is a
non-standard uniserial $R$-module.}

\medskip

\noindent
{\sc Proof}.  Let $E'$ be as above. It may be helpful to regard the
following  as a $\diamondsuit(E')$-like argument carried out in ZFC,
where the $N_\delta$ ($\delta \in E'$)  serve to give the predictions. We will construct the $e^\tau _\sigma\  (\sigma  < \tau  < \delta )$  
satisfying (\dag ) by induction on $\delta $ so that for every $\delta  \in  
E',$ 

\labpar{\#}{there is no sequence of units $\{c_\sigma \colon \sigma  < \delta \} \subseteq  
N_\delta $ such that there exists $c_\delta  \in  R^{*}$ with
$c_\delta  -  e^\delta _\sigma c_\sigma  \in  r_\sigma A\hbox{    for 
all }\ \sigma  < \delta.
$}

Suppose we can do this. We claim that the uniserial module given by the 
$e^\tau _\sigma $ is non-standard. Indeed if not, there exists $\{c_\sigma \colon 
\sigma  < \omega _1\}$ as in (\dag \dag ); then there is a cub $C$ such that 
for all $\delta  \in  C$, $c_\sigma  \in  N_\delta $ if $\sigma  < \delta $; 
but then the construction is contradicted since $C \cap  E'$ is non-empty. 

The definition of the $e^\tau _\sigma $ is routine except for the case when 
$\delta  \in  E'$ and we are defining $e^\delta _\sigma $ assuming the 
$e^\tau _\sigma $ have been defined for $\sigma  < \tau  < \delta $. First 
choose any elements $
\tilde{e}
^\delta _\sigma $ so that  $\{e^\tau_\sigma \colon \sigma < \tau <
\delta\} \cup \{\tilde{e}^\delta_\sigma \colon \sigma < \delta \}$
satisfies (\dag ). If (\#) is true, then we can let 
$e^\delta _\sigma  = 
\tilde{e}
^\delta _\sigma $. Otherwise, we have $\langle d_\sigma \colon \sigma  < 
\delta \rangle  \subseteq  N_\delta $ such that there exists $d_\delta  \in  
R^{*}$ with
$$
d_\delta  -  \tilde{e}^\delta _\sigma d_\sigma  \in  r_\sigma A\hbox{   for 
all }\ \sigma  < \delta .
$$

\noindent
Let $\langle u_\sigma \colon \sigma  < \delta \rangle $ be as in the definition of 
$E'$ (i.e., there is no $f \in  R^{*}$ such that for all  $\sigma  < \delta $, 
$u_\sigma f \in  N_\delta $ mod $r_\sigma A)$. Let $e^\delta _\sigma  = 
u_\sigma 
\tilde{e}
^\delta _\sigma $. We claim that (\#) is satisfied. 

Suppose not; let $\{c_\sigma \colon \sigma  < \delta \} \subseteq  N_\delta $ such 
that there exists $c_\delta  \in  R^{*}$ with
$$
c_\delta  -  e^\delta _\sigma c_\sigma  \in  r_\sigma A\hbox{    for 
all }\ \sigma  < \delta .
$$

\noindent
Let $f = c^{-1}_\delta d_\delta $. Fix $\sigma  < \delta $ and let $\equiv$  
denote congruence mod $r_\sigma A$. Then we have $c_\delta  \equiv  
 e^\delta _\sigma c_\sigma  \equiv   u_\sigma 
\tilde{e}
^\delta _\sigma c_\sigma $ and $d_\delta  \equiv  d_\sigma 
\tilde{e}
^\delta _\sigma $, so $f \equiv  u^{-1}_\sigma c^{-1}_\sigma d_\sigma $. Thus 
$u_\sigma f \equiv  c^{-1}_\sigma d_\sigma  \in  N_\delta $, a contradiction of
the choice of $\langle u_\sigma \colon \sigma  < \delta \rangle $.  \qed

\bigskip

While the condition $\Gamma '_R(J/A) \neq  0$ is thus sufficient for there to be a
non-standard uniserial $R$-module of type $J/A$, we shall see later
(Corollary 15) that it is not necessary.

Let us say that a type $J/A$ is {\it essentially uncountable} if
 for every $\nu  < \omega _1$ there exists $\mu  > \nu $
such that $r_\nu A/r_\mu A$ is uncountable. This will be the case, for example,
if $R$ has uncountable residue field, $J$ is not countably generated and
$A 
= R$.

\bigskip

\noindent
{\bf THEOREM 8.}
$(\aleph _1 < 2^{\aleph _0})$   {\it If $R$ is a valuation domain of 
cardinality $\aleph _1$ and $J/A$ is essentially uncountable, then there is a 
non-standard uniserial $R$-module of type $J/A$.}

\medskip

\noindent
{\sc Proof}.  It suffices to prove that $\Gamma '_R(J/A) \neq  0$; in fact we'll
show that $\Gamma '_R(J/A) = 1$. Without loss of generality we can assume $J = 
\cup _{\nu <\omega _1} r^{-1}_\nu R$ where for all $\nu  < \omega _1$, 
$r_\nu A/r_{\nu +1}A$ has cardinality $\aleph _1$. Let $R = 
\cup _{\alpha <\omega _1}N_\alpha $ be an $\omega _1$-filtration of $R$ by 
subrings. Let $\delta $ be a limit ordinal. Choose a {\it ladder} on
$\delta$, i.e., a strictly increasing 
sequence $\gamma _n\ (n \in  \omega )$ whose limit is $\delta $. We will define 
by induction on $n \in  \omega $ a unit $u_\eta $ for each $\eta  \in  {}^n2$ --- 
the set of all functions from $n = \{0$, 1, \ldots, $n-1\}$ to $2 = \{0,
1\}$ --- such that 
$u_\eta  - u_{\eta \rest k} \in  r_{\gamma _k}A$ for all $k < n$. If $u_\eta $ 
has been defined and $\eta _i \in  {}^{(n+1)}2$ such that $\eta _i\rest n = 
\eta $ and $\eta_i (n) = i\ (i = 0$, 1), then we choose $u_{\eta _i}$ congruent to
$u_\eta $ mod $r_{\gamma _n}A$ and such that $u_{\eta _1}u^{-1}_{\eta _0} 
\notin  N_\delta $ mod $r_{\gamma _{n+1}}A$. This is possible since $N_\delta $
is countable and $r_{\gamma _n}A/r_{\gamma _{n+1}}A$ is uncountable.
(Indeed, choose $ a \in A$ such that $r_{\gamma_n}a \not\equiv y-1$
mod $r_{\gamma_{n+1}}A$ for any $y \in N_{\delta}$; then let
$u_{\eta_0} = u_\eta$ and $u_{\eta_1} = (1 + r_{\gamma_n}a)u_\eta$.)
 Thus for 
each of the $2^{\aleph _0}$ elements $\zeta $ of $^\omega 2$, we have a 
different family $\langle u_{\zeta \rest n}\colon n \in  \omega \rangle $. We claim 
that for at least one  $\zeta $ there is no $f \in  R^{*}$ such that for all
$k$, $u_{\zeta \rest k}f \in  N_\delta $ mod $r_{\gamma _k}A$ ---
which will show $\delta \notin E'$. Indeed, if there
were no such $\zeta $, then since $R^{*}$ has cardinality $\aleph _1
< 2^{\aleph_0}$, there would
be $\gz_0 \neq \gz_1$ and  $f$  satisfying

$$
u_{\zeta _j\rest k}f \in  N_\delta \hbox{ mod }\ r_{\gamma _k}A
$$

\noindent
for all $k \in \go$ and $j \in \{0, 1\}$. Let $n$ be such that $\zeta _0\rest n =
\zeta _1\rest n$ but $\zeta _0(n) \neq  \zeta _1(n)$.  It is then easy to obtain a contradiction of our construction
of $u_{\eta _i} \  (i = 0,1)$ if we let let $\eta  = 
\zeta _0\rest n$ and $\eta_i = \zeta_i \rest (n+1)$.  \qed

\bigskip

\noindent
{\bf COROLLARY 9}. {\it Let $R$ be a valuation domain of cardinality
$\hao$. If $J/A$ is essentially uncountable and $\Gamma _R(J/A) = 1$,
then there is a non-standard uniserial $R$-module of type $J/A.$}

\medskip

\noindent
 {\sc Proof}.  If CH holds, this is by Proposition 3; if CH fails this is
by Theorem 8.  \qed

\bigskip

We will deal with the essentially countable case in a later paper; in
fact, Corollary 9 can be proved without the hypothesis of essential
uncountability, but this requires a different construction in the
essentially countable case.
This means that Osofsky's original conjecture is vindicated for
``natural'' valuation domains, where $\Gamma_R(J/A)$ has to be either $0$ or
$1$.

\bigskip

\bigskip

\noindent
{\bf A construction of non-standard uniserials.} 

\medskip

As an application of Theorem 7, we will show how to construct
non-standard uniserials (of type $Q/R)$ in ZFC over  certain
concretely given valuation domains
R. This may be regarded as a generalization of [O2] in that we define
such  $R$ with any given non-zero value for
$\Gamma _R(Q/R)$. We borrow an idea from [O2]
 (the definition of $u_\sigma )$  in order to show that $\Gamma'_R(Q/R) >
0.$
First we define the ring.

\bigskip

\noindent
{\bf Definition 10.}   Let $G$ be the ordered abelian group which is the direct
sum $\oplus _{\alpha <\omega _1}{\open Z}\alpha $ ordered 
anti-lexicographically; that is, $\Sigma _\alpha  n_\alpha \alpha  > 0$ if and 
only if $n_\beta  > 0$, where $\beta $ is maximal such that $n_\beta  \neq  0$.
(This is the group used in [FS] and [O2] and denoted $\Gamma (\omega
_1)$ in [BS1; {\S}3]). In particular, the
basis elements have their natural order and if $\alpha  < \beta $, then 
$k\alpha  < \beta $ in $G$ for any $k \in  {\open Z}.)$  Let $G^+ = \{g \in  G\colon
g \geq  0\}$ and for each $\gamma  < \omega _1$, let $G^+_\gamma $ be the 
submonoid $\{g \in  G\colon 0 \leq  g < \gamma \}$ of $G^+.$ 

Fix a field $K$ of cardinality $\leq  \aleph _1$ and let 

$$\hat{R}  = \{\Sigma _{g\in \Delta } k_gX^g\colon k_g \in  K, \Delta \hbox{ a
well-ordered subset of } G^+\},$$

\noindent
(= $K[[G]]$; cf. [FS; pp. 24f]) 

$$\hat{R} _{\alpha +1} = \{\Sigma _{g\in \Delta } k_gX^g \in  \hat{R} \colon
\Delta \subseteq  G^+_\alpha  \}$$

\noindent
and for $\gamma $ a limit ordinal,

$$\hat{R} _\gamma  = \cup _{\alpha <\gamma } \hat{R} _\alpha. $$

Given an element $y = \Sigma _{g\in \Delta } k_gX^g$ of $\hat{R} $, let
$\supp(y) = \{g \in  \Delta \colon k_g \neq  0\}$; let $\psupp(y) = \{\alpha  \in  
\omega _1\colon \exists g \in  \supp(y)$ whose projection on ${\open Z}\alpha $ is 
non-zero$\}$. If $X \subseteq  G$, then $y\rest X$ is defined to be 
$\Sigma _{g\in X\cap \Delta } k_gX^g$. Let $y\rest \ga  = y\rest \{g \in  G\colon g
< \ga \}$, the {\em truncation} of $y$ to $\ga$. Note that $y \in
\hat{R}_{\ga + 1}$ if and only if $y = y\rest\ga$.

For $Y \subseteq  \hat{R} $, let $\langle Y\rangle _{val}$ denote the 
valuation subring generated by $Y$, i.e., the intersection with $\hat{R} $ 
of the quotient field of the subring generated by $Y.$ 

For any stationary subset $S$ of $\omega _1$, let $R_S$ be the subring of 
$\hat{R} $ consisting of all elements $y$  of $\hat{R} $ such that 
$\supp(y)$ is countable and no member of $S$ is a limit point of $\psupp(y)$ 
(in the order topology).

Then $R_S$ is a valuation domain since 
$ \psupp(xy^{-1}) \subseteq \psupp(x)\ \cup\ \psupp(y)$. The cardinality
of $R_S$ is $\leq 2^{\aleph_0}$.

\bigskip

\noindent
{\bf THEOREM 11.} (CH)  {\it For every field $K$ of cardinality $\aleph _1$ and
every stationary subset $S$ of $\omega_1$, there is a valuation domain $R$  
of cardinality $\hao$
 with residue field $K$ and quotient field $Q$
such that $\Gamma _{R}(Q/R) = \tilde{S}$ and there
is a non-standard uniserial $R$-module of type $Q/R.$}

\medskip

\noindent
{\sc Proof}. Let $R$ be $R_S$ as defined above. Let $r_\sigma  = X^{\sigma}$, so that $Q =
\cup _{\sigma <\omega _1} r^{-1}_\sigma R$. First notice that $\Gamma
_{R}(Q/R) \geq \tilde{S}$
since for any $\delta \in S \cap \lim(\omega_1)$, if $\gamma_n \ (n
\in \omega)$ is a ladder on $\delta$ and $y_n$ is defined to be

$$\Sigma _{i = 0} ^{n-1} X^{\gamma_i}$$

\noindent
then $\{y_n \colon n < \omega \}$ is a Cauchy sequence in $R/(R :
J_\delta)$ which does not have a limit, since a putative limit has
$\delta$ as  a limit point of its p-support. Moreover, $\Gamma _{R}(Q/R)
\leq \tilde{S}$
since for $\delta \notin S$,  a Cauchy
sequence in  $R/(R : J_\delta)$ {\it will} have a limit in $R$
because 
$\delta$ is allowed to be the limit point of the p-support of an
element of $R$.

To show that there is a non-standard uniserial of type $Q/R$, it suffices
by Theorem 7 to
show that $\Gamma '_R(Q/R) \geq  \tilde{S}$. Let us fix an
$\omega _1$-filtration of $R$ by subrings, $\{N_\alpha \colon \alpha  \in  
\omega _1\}$. Let $\delta  \in  S$ and let $\langle \gamma _n\colon n \in  
\omega \rangle $ be a ladder on $\delta $.
Let $L$ be the subfield of $K$ generated by the coefficients of the members of 
$N_\delta $. Then $L$ is a countable subfield of $K$ so we can
inductively choose $t_n \in  K 
\setminus  L(t_0, \ldots, t_{n-1})$. Let
$$
u_\sigma  = 1 + \Sigma ^{n-1}_{i=0} t_iX^{\gamma _i}
$$

\noindent
for all $\sigma $ with $\gamma _{n-1} <  \sigma  \leq \gamma _{n}$. Then clearly 
$\langle u_\sigma \colon \sigma  < \delta \rangle  \in  {\cal T}^\delta _{Q/R}$.  It
suffices to show that there is no $f \in  R^{*}$ such that for all $\sigma  < 
\delta $, $u_\sigma f \in  N_\delta $ mod $X^\sigma R$.
Suppose there is such an $f$;
 by definition of $R_S$ and since $\delta  
\in  S$, there exists $n$ such that $\sup (\psupp(f \rest \delta)) \leq \gamma
_{n - 2}$. Say $f \rest \delta = \Sigma_{g \leq \gamma_{n-2}} k_g X^g$.
But then, because of the properties of the value group --- in
particular because $2\gamma_{n-2} < \gamma_{n-1}$ --- we
see that
the coefficient of $X^{\gamma _{n-1}}$ in $u_{\gamma _n}f$ is
$t_{n-1}k_0$; similarly
 the coefficient of $X^{\gamma _{n}}$ in $u_{\gamma _{n+1}}f$
is $t_n k_0$.
By assumption, these coefficients belong to $L$, and we easily obtain a contradiction of the choice of the $t_i$.
  \qed 

\medskip

\noindent
{\bf Remarks.} (1) The hypothesis of CH was used only to insure that
$R_S$ has cardinality $\hao$. The ring defined by Osofsky in [O2] has
cardinality $\hao$ without invoking CH, and the proof above applies to
it. (It satisfies $\Gamma'(Q/R) = 1$.) Bazzoni has pointed out that the
proof of Theorem 7 can be adapted to prove Theorem 11 for the rings
$R_S$, even when $R_S$ has cardinality $2^{\ha_0} >
\hao$ and so cannot have an $\goo$-filtration by countable subrings:
choose an $\goo$-filtration of $K$ by countable subfields $L_\ga$ and
define $N_\ga$ to be the subrings of elements whose coefficients come
from $L_\ga$.

(2) The main result of [O1] is very general in that it deals with
arbitrary types, not just $Q/R$. We have not been able to duplicate this
generality, but the method of proof of Theorem 11 can be used to prove,
in ZFC, that there exist non-standard
uniserials of other types than $Q/R$, e.g., ones like those given in Examples 1 -
6 of [BS1; pp. 301--305].

\bigskip

\bigskip

\noindent
{\bf Axiom ($S$).} 

\medskip

As we have observed, Lemma 5 says that the condition $\Gamma _R(J/A) \neq  0$ 
is necessary for the existence of a non-standard uniserial of type
$J/A$. 
 We
now aim to prove that the condition $\Gamma _R(J/A) \neq  0$ is {\it not} 
provably sufficient if we assume only GCH (rather than V = L), i.e., the ``if"
direction of the conjecture fails in some model of ZFC + GCH. By
Proposition 3 and Theorem 7,
we will need to consider $R$ and $J/A$ where 
$\Gamma' _R(J/A) = 0$ and 
$\Gamma _R(J/A) \neq 
 1$. Analogies with the Whitehead Problem suggest the use of the 
following principle, known as Ax($S$), which has been shown by the second 
author to be consistent with ZFC + GCH (and with ZFC + $\neg$ CH):  
\begin{quotation}
\noindent
there is a stationary and co-stationary subset $S$ of $\omega _1$ such that for
all proper posets ${\open P}$ of cardinality $\aleph _1$ which are $(\omega _1 
\setminus  S)$-complete, and for all families ${\cal D}$ of $\aleph _1$ dense 
subsets of ${\open P}$, there is a ${\cal D}$-generic subset of ${\open P}.$
\end{quotation}
\noindent
(See [EM; pp. 170-173] for the necessary definitions.) 

\smallskip

 We will be able to prove that in a model of Ax($S$) + CH, 
 there exists a valuation domain $R$ of cardinality $\aleph _1$ with 
$\Gamma _R(Q/R) = \tilde{S}$ such that every uniserial of type $Q/R$ is standard.

\bigskip

\noindent
{\bf Definition 12.}   For any countable subring $N$ of a valuation domain $R$ 
of cardinality $\aleph _1$, let 

$${\cal T}^\delta _{J/A}(N) = \{\langle \bar u_\sigma \colon \sigma  < \delta \rangle
\colon \forall \sigma < \gt  < \delta (u_\sigma  \in  N \hbox{ and }
u_\gt - u_\gs \in r_\gs A)\}$$
where $\bar u_\gs$ denotes the coset $u_\gs + r_\gs A$. Thus
${\cal T}^\delta _{J/A}(N)$ consists of sequences
$\langle \bar u_\sigma \colon \sigma  < \delta \rangle$ where
$\langle  u_\sigma \colon \sigma  < \delta \rangle \in {\cal
T}^\gd_{J/A}$ and all the $u_\gs$ belong to $N$. Usually we shall abuse
notation and write $\langle  u_\sigma \colon \sigma  < \delta \rangle$
for an element of ${\cal T}^\delta _{J/A}(N)$; but notice that
$\bar u_\gt$ determines $\bar u_\gs$ for all $\gs < \gt$. Let

$${\cal L}^\delta _{J/A}(N)= \{\langle \bar u_\sigma \colon
 \sigma  < \delta \rangle \in
 {\cal T}^\delta _{J/A}(N) \colon
\langle  u_\sigma \colon
 \sigma  < \delta \rangle \in
{\cal L}^\delta _{J/A}\}.$$

\noindent
${\cal T}^\delta _{J/A}(N)$ is given the tree topology, i.e., if $u = 
\langle \bar u_\sigma \colon \sigma  < \delta \rangle  \in  {\cal T}^\delta _{J/A}(N)$, a
basis of neighborhoods of $u$ is the family $\{[u]_\nu \colon \nu  < \delta \}$, 
where $[u]_\nu  = \{\langle \bar v_\sigma \colon \sigma  < \delta \rangle  \in
{\cal T}^\delta _{J/A}(N)\colon \bar v_\sigma  = \bar u_\sigma $ for all $\sigma  \leq
\nu \}.$

\bigskip

>From now on $S$ always denotes the stationary and co-stationary subset of 
$\omega _1$ which is asserted to exist by Ax($S$).

\bigskip

\noindent
{\bf THEOREM 13.}  (Ax($S$))   {\it Suppose that $R$ is a valuation domain of 
cardinality $\aleph _1$ such that $\Gamma _R(J/A) \subseteq  \tilde{S}$ and 
$\Gamma '_R(J/A) = 0$. Suppose in addition that for some $\omega _1$-filtration
of $R$ by subrings, $R = \cup _{\alpha \in \omega _1} N_\alpha $, the following
holds for all $\delta  \in  S$ :}
\labpar{\&}{
  {\it for every open subset $U$ of ${\cal T}^\delta _{J/A}(N_\delta )$, $U
\cap  {\cal L}^\delta _{J/A}(N_\delta )$ is a non-meager (i.e., second category) 
subset of ${\cal T}^\delta _{J/A}(N_\delta ).$}}

\noindent
{\it Then every uniserial $R$-module of type $J/A$ is standard.}

\medskip

\noindent
{\sc Proof}.  Let $U$ be a uniserial $R$-module of type $J/A$ given by a family 
$\{e^\tau _\sigma \colon \sigma  < \tau  < \omega _1\}$ as in (\dag ). Let the poset
${\open P}$ consist of all sequences $p = \langle \bar p_\sigma \colon \sigma  \leq
\mu \rangle $ where $\mu $ is a countable ordinal (denoted $\ell (p)$ and 
called the {\it length} of $p)$ such that  $$\forall \sigma  < \tau  \leq
\mu (p_\sigma  \in  R^{*}, \hbox{ and } p_\tau  - e^\tau _\sigma p_\sigma  \in
r_\sigma A).$$
The partial ordering is the natural one of extension of
sequences. In
an abuse of notation we shall write the elements of ${\open P}$
as if they were sequences of elements of $R$ rather than sequences of cosets.
Note that ${\open P}$ is a tree and has the following properties:

\begin{quotation}

(i) if $p$, $q \in  {\open P}$, $p_\tau  = q_\tau $ for some $\tau $ and $c$ is
a sequence of length $= \ell (q) \geq  \tau $, such that $c_\sigma  \equiv  
p_\sigma $ (mod $r_\sigma A)$ for all $\sigma  \leq  \tau $, and $c_\sigma  
\equiv  q_\sigma $ (mod $r_\sigma A)$ for all $\tau  \leq  \sigma  \leq  
\ell (q)$, then $c \in  {\open P};$ 

\smallskip

(ii) if $p \in  {\open P}$ and $e \in  R^{*}$, then $\langle ep_\sigma \colon \sigma  
\leq  \ell (p)\rangle  \in  {\open P}.$ 
\end{quotation}

For each $\nu  < \omega _1$, let $D_\nu  = \{p \in  {\open P}\colon p$ has length
$\geq  \nu \}$. Then $D_\nu $ is dense in ${\open P}$ because given $q = 
\langle q_\sigma \colon \sigma  \leq  \mu \rangle  \in  {\open P}$ (where we can 
assume $\mu  < \nu )$, define $p_\sigma  = q_\sigma $ if $\sigma  \leq  \mu $ 
and $p_\sigma  =(e^\nu _\sigma )^{-1}e^\nu _\mu q_\mu $ if $\mu  \leq  \sigma  
\leq  \nu $. Then if $\sigma  < \tau  \leq  \nu $, one may easily check that 
modulo $r_\sigma A,$
$$
p_\tau  - e^\tau _\sigma p_\sigma  \equiv  ((e^\nu _\tau )^{-1} - 
e^\tau _\sigma(e^\nu _\sigma )^{-1} )e^\nu _\mu q_\mu  \equiv  
0\cdot e^\nu _\mu q_\mu 
$$

\noindent
if $\mu  \leq  \sigma  < \tau  \leq  \nu $, and similarly if $\sigma  < \mu  < 
\tau  \leq  \nu $. Thus $\langle p_\sigma \colon \sigma  \leq  \nu \rangle $ is an 
element of ${\open P}$ extending $q$, so $D_\nu $ is dense in ${\open P}.$ 

It suffices then to prove that there is a $\{D_\nu \colon \nu  < 
\omega _1\}$-generic subset of ${\open P}$; and for this, given the 
set-theoretic hypothesis, it suffices to prove that ${\open P}$ is an 
$(\omega _1 \setminus  S)$-complete proper poset.

Let $\kappa $ be large enough for ${\open P}$; let ${\cal C}$ be the cub of all
countable elementary submodels $N$ of $(H(\kappa )$, $\in$, $\omega _1$,
$\{e^\tau _\sigma \colon \sigma  < \tau  < \omega _1\}$, $R$, ...) such that $N =
\cup _{i\in \omega }N_i$ where $N_i \prec  N_{i+1}$ and $N_i \cap  \omega _1 <
N_{i+1} \cap  \omega _1$. Fix $N \in  {\cal C}$ and let $\delta  = N \cap
\omega _1$, $\delta _i = N_i \cap  \omega _1$. To
show that ${\open P}$ is $(\omega _1 
\setminus  S)$-complete consider a chain
$$
p^0 \leq  p^1 \leq  \ldots\hbox{.}
$$

\noindent
of elements $p^n = \langle p^n_\sigma \colon \sigma  \leq  \ell (p^n)\rangle $ of 
${\open P} \cap  N$  such that for 
each dense subset $D$ of ${\open P}$ which belongs to $N$ there exists $n$ with
$p^n \in  D$.
For all $n$ there exists $m_n$ such that  $\ell (p^{m_n}) \geq  \delta _n$ (because
$D_{\delta_n} \in N_{n+1}$). For each
$\sigma < \delta$, let  $p_\sigma = p^m_\sigma $ for any $m$ such
that $\ell(p^m) \geq \sigma$. 
 Then $\{e^\delta _\sigma p_\sigma \colon \sigma  < \delta \}$ is a 
Cauchy sequence in $R/(A : J_\delta )$, so if $\delta  \notin  S$, it has a 
limit $p_\delta  \in  R$, and then $\langle p_\sigma \colon \sigma  \leq  
\delta \rangle $ is a member of ${\open P}$ which extends each $p^n$.

So far, we have used just the hypothesis that $\Gamma _R(J/A) \leq  \tilde{S}$.
To show that ${\open P}$ is proper, we need to use the additional hypotheses on
R. Since the intersection of two cubs is a cub, we can 
assume, without loss of generality, that $N  \cap R = N_\delta $ (the $\delta$th
 member of
the $\omega _1$ filtration of $R$ given in the statement of the theorem).

Let $q \in  {\open P} \cap  N$; we must find $p \geq  q$ such that $p$ is 
$N$-generic. This is no problem if $\delta  \notin  S$ since ${\open P}$ is 
$(\omega  \setminus  S)$-complete, so assume $\delta  \in  S$. Choose $d = 
\langle d_\sigma \colon \sigma  \leq  \delta \rangle  \in  {\open P}$ of length 
$\delta $ and let $u_\sigma  = e^\delta _0(e^\sigma _0)^{-1}d_\sigma $.
Then $\langle u_\sigma \colon \sigma < \delta \rangle \in {\cal T}^\delta_{J/A}$.
 Since 
$\Gamma '_R(J/A) = 0$ we can assume that there exists $\langle
f_\sigma \colon \sigma < \delta \rangle \in  
{\cal L}^\delta _{J/A}$ with limit $f_\delta$  so that $u_\sigma f_\sigma  \in  N_\gd$ for all $\sigma  <
\delta $. Then, replacing $d$ by $\langle e^\delta _0d_\sigma f_\sigma \colon 
\sigma  \leq  \delta \rangle $, we can assume that $d_\sigma $ belongs to $N_\gd$
for all $\sigma  < \delta $ (since $e^\delta _0d_\sigma f_\sigma  = 
e^\sigma _0u_\sigma f_\sigma $; note that $e^\sigma _0 \in  N_\gd$  for all $\sigma
< \delta .)$ Also, by (i) and (ii), we can assume that $d_\sigma  = q_\sigma $ 
for all $\sigma  \leq  \ell (q)$. 
(Note that $\ell(q) < \delta$ since $q \in N$ and $N \cap \omega_1 =
 \delta$.)

Let $U$ be the open subset
$$
\{\langle u_\sigma \colon \sigma < \delta \rangle \in  {\cal
T}^\delta _{J/A}(N_\gd) \colon u_\sigma  = 1\hbox{ for }\ \sigma  \leq
\ell (q)\}\hbox{.}
$$

\noindent
Given an element $u = \langle u_\sigma \colon \sigma  < \delta \rangle $ of 
${\cal T}^\delta _{J/A}(N_\gd)$, we will let $du$ denote $\langle d_\sigma u_\sigma \colon
\sigma  < \delta \rangle $; notice that for all $\sigma  < \tau  < \delta $ we 
have
$$
d_\tau u_\tau  - e^\tau _\sigma d_\sigma u_\sigma  = d_\tau u_\tau  - 
d_\tau u_\sigma  + d_\tau u_\sigma  - e^\tau _\sigma d_\sigma u_\sigma  \in  
r_\sigma A,
$$

\noindent
so $du$ belongs to ${\open P}$. For each dense subset $D$ of ${\open P}$ which belongs to $N$, let
$$
D' = \{u \in  {\cal T}^\delta _{J/A}(N_\gd) \colon du \hbox{ is not an extension of a member
of }\ D \cap  N\}.
$$

\noindent
Then $D'$ is a nowhere dense subset of ${\cal T}^\delta _{J/A}(N_\gd)$ since for any
basic open subset $[u]_\nu $, we can choose $u' = \langle u'_\sigma \colon \sigma  
\leq  \ell (u')\rangle  \geq  \langle d_\sigma u_\sigma \colon \sigma  \leq  
\nu \rangle $ which belongs to $D \cap  N$, and then 
$\langle d^{-1}_\sigma u'_\sigma \colon \sigma  \leq  \ell (u')\rangle $ determines 
a basic non-empty open subset of $[u]_\nu $ which is disjoint from $D'$. 
(We are using here the fact that in the elementary submodel $N$, $D \cap N$ is
a dense subset of ${\open P} \cap N$; moreover note that $\langle
\overline{d_\sigma u_\sigma} \colon \sigma  \leq
\nu \rangle $ belongs to ${\open P} \cap N$.)
Hence 
by (\&), there is an element $w = \langle w_\sigma \colon \sigma  < \delta \rangle $
of $U \cap  {\cal L}^\delta _{J/A}$ which does not belong to $D'$ for any $D 
\in  {\open P} \cap  N$. Let $w_\delta $ be such that $w_\delta  - w_\sigma  
\in  r_\sigma A$ for all $\sigma  < \delta $. Then $\langle d_\sigma w_\sigma \colon
\sigma  \leq  \delta \rangle $ is the desired $N$-generic element extending 
$q$.  \qed

\bigskip

Now the following result shows that the conjecture (\$) fails in a model of 
Ax($S$) + CH. The hypothesis of CH is essential: see Theorem 8.

\bigskip

\noindent
{\bf THEOREM 14.}  (Ax($S$) + CH)   {\it For any field $K$ of cardinality 
$\leq  \aleph _1$, there is a valuation domain $R$ of cardinality $\aleph _1$ 
with residue field $K$ and a stationary, co-stationary subset $S$ of 
$\omega _1$ such that $\Gamma _R(Q/R) = \tilde{S}$, but every uniserial $R$-module of
type $Q/R$ is standard.}

\medskip

\noindent
{\sc Proof}.  Let $\hat{R} $ be as in Definition 10; it suffices to show there is a
valuation subring $R$ of $\hat{R} $ of cardinality $\aleph _1$ which 
satisfies the hypotheses of Theorem 13 for the type $Q/R$. We let $r_\nu  = 
X^\nu $ which will be in R. (Thus $Q$ will be $\cup
_{\nu<\omega _1}r^{-1}_\nu R.)$  For each $\gamma  \in  \omega _1$, let 
$W_\gamma  = \{y \in  \hat{R} _\gamma \colon \supp(y)$ is finite and the 
coefficients of $y$ belong to the prime subfield of $K\}.$ 

We will now define by induction on $\alpha  < \omega _1$ a continuous chain 
$\{Y_\alpha \colon <  \omega _1\}$ of valuation subrings of $\hat{R} $ such that 
for all $\alpha$:  $W_\alpha  \subseteq  Y_\alpha  \subseteq
\hat{R} _\alpha $; $Y_\alpha  \cap  R^{*} = Y^*_{\alpha}
$. Moreover, the $Y_\ga$ will be closed under truncation,
that is, if $y \in Y_\ga$ and $\gb < \ga$, then $y\rest\gb \in
Y_{\gb + 1}$.
We will then define $R$ to be the union of the $Y_\alpha $.

First of all notice that  each element of ${\cal T}^\gd_{Q/R}(N)$
is uniquely represented by a sequence $u =
\langle u_\sigma \colon \sigma < \delta \rangle$ where
for all $\gs < \gd$,
$u_\gs \in N$,  and $u_\gt\rest\gs = u_\gs$
whenever $\gt \geq \gs$. We shall consider these to be the elements of
${\cal T}^\gd_{Q/R}(N)$. Moreover, given $u =  \langle u_\sigma \colon
\sigma < \delta \rangle \in {\cal T}^\gd_{Q/R}(Y_\gd)$,
there is a unique element of $\hat{R} _{\delta +1}$, which we will denote 
$u_\delta $, which represents the limit of $ \langle u_\sigma \colon \sigma  <
\delta \rangle $ in $\hat{R}/\cap_{\nu < \gd} r_\nu\hat{R}$. If $u_\gd \notin Y_{\gd + 1}$,
then, by closure under truncation, the sequence $ \langle u_\sigma \colon \sigma  <
\delta \rangle $ will not have a limit in $R$.
Given a subset $Z$ of 
${\cal T}^\delta _{Q/R}(Y_\delta )$, we will denote by $Z_{[\delta ]}$ the set 
$\{u_\delta \colon u \in  Z\} \subseteq  \hat{R} _{\delta +1}.$

Let $Y_0 = \langle \{1\}\rangle _{val}$. Suppose that $Y_\beta $ has been 
defined for all $\beta  < \gamma $. If $\gamma $ is a limit ordinal, let 
$Y_\gamma  = \cup _{\beta <\gamma } Y_\beta $.  If $\gamma  = \delta  + 1$ and 
$\delta $ is a successor ordinal, let $Y_\gamma $ be the valuation subring of 
$\hat{R} $ generated by $Y_\delta  \cup  W_{\delta +1}$. If $
\delta  \notin  S$, let $Y_\gamma $ be the valuation subring generated by 
$Y_\delta $ together with the $u_\delta $ for each $u \in  
{\cal T}^\delta _{Q/R}(Y_\delta )$.  

If $\delta $ is a limit ordinal  in $S$, we must satisfy the 
hypotheses of Theorem 13 and make sure $R/\cap_{\nu < \gd} r_\nu R$ is not complete. Let
$\{(N_\beta $, $\Phi _\beta $, $V_\beta )\colon \beta  < 2^{\aleph _0} = 
\aleph _1\}$ be an enumeration of all triples such that $N_\beta $ is a 
countable subring of $Y_\delta $, $V_\beta $ is a basic open  set in
${\cal T}^\delta _{Q/R}(N_\beta )$, and $\Phi _\beta $ is the complement of a 
countable union of closed nowhere dense subsets of 
${\cal T}^\delta _{Q/R}(N_\beta )$. (This enumeration is possible since each 
${\cal T}^\delta _{Q/R}(N_\beta )$ has a countable base of open sets.) Also 
enumerate ${\cal T}^\delta _{Q/R}(Y_\delta )$ as $\{u^\alpha =
\langle u^\alpha_\sigma \colon \sigma < \delta \rangle \colon \alpha  < 
2^{\aleph _0}\}$. We will  inductively define  elements $y_\alpha $ of
$\hat{R} _{\delta +1}\ (\alpha  < 2^{\aleph _0})$ and then let $Y_{\delta +1} 
= \langle Y_\delta  \cup  \{y_\alpha \colon \alpha  < 
2^{\aleph _0}\}\rangle _{val}$. Fix an element $c = \langle c_\sigma
\colon \sigma < \delta \rangle$ of 
${\cal T}^\delta _{Q/R}(Y_\delta )$ and let $c_\delta \in
\hat{R}_{\delta + 1}$ be its limit, as above. We also require that
$c$ is chosen so that $c_\delta \notin \hat{R}_\gd$.
 The elements 
$y_\alpha $ will be chosen with the following properties for all $\alpha  < 
2^{\aleph _0}$: 

\begin{quotation}

(i) if $\alpha  = 2\beta  + 1$, then $y_\alpha  = u_\delta $ for some $u \in  
\Phi _\beta  \cap  V_\beta ;$ 

(ii) if $\alpha  = 2\beta $, then $y_\alpha  = f_\delta $ for some $f =
\langle f_\gs \colon \gs < \gd \rangle \in
{\cal T}^\delta _{Q/R}(Y_\delta )$ such that $u^\beta _\sigma f_\sigma  \in  
W_\delta $ for all $\sigma  < \delta ;$ 

(iii) $c_\delta \notin  \langle Y_\delta  \cup  \{y_\grg \colon \grg  <
\alpha \}\rangle _{val}.$
\end{quotation}

\noindent
Then (i) (resp. (ii)) will imply that (\&) holds (resp. $\Gamma '_R(Q/R) = 0)$,
and (iii) will insure that $R/\cap_{\nu < \gd} r_\nu R$ is not complete.

Suppose that $y_\gamma $ has been defined for all $\gamma  < \alpha $ and that 
$\alpha  = 2\beta +1$. Temporarily let $R' = \langle Y_\delta  \cup  
\{y_\gamma \colon \gamma  < \alpha \}\rangle _{val}$. We claim that: 
\labpar{\hbox{14.1}}{
 there is a subset $Z$ of $\Phi _\beta  \cap  V_\beta $ of cardinality 
$> 1$ such that $Z_{[\delta ]}$ is algebraically independent over $R'$.}

\noindent
If so,  there exists $u \in  Z$ such that $c_\delta \notin  \langle R' \cup  
\{u_\delta \}\rangle _{val}$  and then we can let $y_\alpha  = u_\delta $.
 For otherwise, we have $u \neq  u'$ in $Z$ and non-zero
polynomials $f_1(T)$, $f_2(T)$, $g_1(T)$, $g_2(T) \in  R'[T]$ such that $c_\delta= 
f_1(u_\delta )/g_1(u_\delta )$ and $c_\delta = f_2(u'_\delta )/g_2(u'_\delta )$. But
then the polynomial $f_1(T_1)g_2(T_2) - f_2(T_2)g_1(T_1) \in  R'[T_1$,\ $T_2]$ 
is non-zero (because $c_\gd \notin R'$) and
shows that $\{u_\delta, u'_\delta \}$ --- and hence $Z_{[\delta ]}$ --- is not 
algebraically independent over $R'$.  
To prove (14.1), we show 
\labpar{\hbox{14.2}}{there is a subset $Z'$ of $\Phi _\beta  \cap  V_\beta $ of cardinality 
$2^{\aleph _0}$ such that $Z'_{[\delta ]}$ is algebraically independent over 
$Y_\delta $.}

\noindent
This suffices, by the additivity of transcendence degree, since $|\{y_\gamma \colon 
\gamma  < \alpha \}| < 2^{\aleph _0}.$ 

In order to prove (14.2), we write the complement of $\Phi _\beta $ as the 
union,  $\cup _{n \in \omega} C_n$, of countably many closed nowhere dense sets in 
${\cal T}^\delta _{Q/R}(N_\beta )$. Choose a ladder
$\{\tau _n\colon n \in  \omega \}$ on $\delta $ and inductively define a 
tree of basic clopen subsets $U_\eta  \subseteq  V_\beta\ (\eta  \in  {}^n2)$ so 
that $U_\eta  \subseteq  U_\zeta $ if $\zeta  = \eta \rest \ell
(\zeta )$ and 
$U_\eta  \subseteq$  the complement of $C_{\ell (\eta )}$. Moreover, we choose 
the $U_\eta $ so that for each $n$, if $^n2 = \{\eta _1$, \ldots, $\eta _{2^n}\}$, then there 
are ordinals $\gamma _i$ with $\tau _n < \gamma _0 < ..$. $< \gamma _{2^n}$
such that for every $z_i \in  U_{\eta _i}\ (i = 1$, \ldots, $2^n)$ and $g \in  
G^+$, if $\gamma _1 \leq  g \leq  \gamma _{2^n}$, then $g \in  \supp(z_i)$ if 
and only if $g = \gamma _{i-1}$. Then one may prove by induction on $m$
 that if 
$f(T_1$, \ldots, $T_m) \in  Y_\delta [T_1$, \ldots, $T_m]$ is a non-zero polynomial 
whose coefficients are in $\hat{R}_{\tau _n}$ and if $z_j \in  U_{\eta _j}$, then
$f(z_1$, \ldots, $z_m) \in \hat{R}_{\gamma_{m+1}} \setminus \{ 0\}$. Hence if we choose for each $\xi  \in
{}^\omega 2$, $z_\xi  \in  \cap _{n \in \omega} U_{\xi \rest n}$, then $Z_{[\delta ]} = 
\{z_\xi \colon \xi  \in  {}^\omega 2\}$ is algebraically independent over $Y_\delta .$

Finally consider the case when $\alpha  = 2\beta $. It suffices to show that 
there is a subset $Z$ of ${\cal T}^\delta _{Q/R}(Y_\delta )$ of cardinality 
$2^{\aleph _0}$ so that $Z_{[\delta ]}$ is algebraically independent over 
$Y_\delta $ and  each $f \in  Z$ has the form $\langle  
w_\sigma(u^\beta _\sigma )^{-1} \colon \sigma < \delta \rangle$
where 
$\langle w_\sigma \colon \sigma <  \delta \rangle \in {\cal T}^\delta _{Q/R}(Y_\delta )$ and for each $\sigma  < \delta $, $w_\sigma 
\in  W_\delta $.  Now by a
tree argument similar to that above,
 there is a subset $Z'$ of ${\cal T}^\delta _{Q/R}(Y_\delta )$ of 
cardinality $2^{\aleph _0}$ consisting of elements
$w = \langle w_\gs \colon \gs < \gd \rangle$ such that for all $\gs$ the
 coefficients of $w_\gs$ belong to
$\{0$, $1\}$ and such that $Z'_{[\delta ]}$ is algebraically independent over 
$Y_\delta $. Then by an argument on transcendence degree, there is a subset $S$
of $\{w_\delta (u^\beta  )^{-1}_\gd \colon w \in  Z'\}$ of cardinality
$2^{\aleph _0}$ which is algebraically independent over $Y_\delta $. So let
$$
Z = \{ \langle w_\sigma (u^\beta _\sigma )^{-1}\colon\ \sigma  <
\delta \rangle \colon w_\delta (u^\beta )^{-1}_\gd
\in  S\}\hbox{.  \qed }
$$

\bigskip

The following corollary shows that the converse of Theorem 7 is not
provable in ZFC + GCH.

\bigskip

\noindent
{\bf COROLLARY 15.}  (V = L)   {\it There is a valuation domain $R$
of cardinality $\hao$ such that
$\Gamma '_R(Q/R) = 0$ and there is a non-standard uniserial $R$-module of type 
$Q/R.$}

\medskip

\noindent
{\sc Proof}.  Let $R$ be the valuation domain constructed as in the proof of Theorem 
14 for some stationary  S.
(The construction requires only CH.) Then $\Gamma '_R(Q/R) = 0$. Also 
$\Gamma _R(Q/R) = \tilde{S} \neq  0$, so by Proposition 3, there is a non-standard
uniserial $R$-module, since  V = L, and hence $\diamondsuit _{\omega _1}(S)$ 
holds.  \qed

\bigskip

\bigskip

\noindent
{\bf Valuation domains of cardinality $> \aleph _1.$} 

\medskip

Finally we want to observe that the hypothesis that $R$ has cardinality 
$\aleph _1$ is essential for Proposition 3; that is, it is not enough to 
require that $J$ is generated by $\aleph _1$ elements (in which case the 
definition of $\Gamma _R(J/A)$ makes sense). In fact, we have the following 
theorem of ZFC.

\bigskip

\noindent
{\bf THEOREM 16.}   {\it There is a valuation domain $R$ of cardinality 
$2^{\aleph _1}$ such that $Q$ is generated (as $R$-module) by $\aleph _1$ 
elements and $\Gamma _R(Q/R) = 1$, but every uniserial $R$-module of type $Q/R$
is standard.}

\medskip

\noindent
{\sc Proof}.  Let $G$ be the free abelian group on $\omega _1 \times  2^{\aleph _1}$.
We order $\omega _1 \times  2^{\aleph _1}$ lexicographically and then order $G$
anti-lexicographically (cf. Definition 10). Thus in $G$, $(\alpha $, $i) < 
(\beta $, $j)$ if and only if $\alpha  < \beta $ (in $\omega _1)$ or $\alpha  =
\beta $ and $i < j$ (in $2^{\aleph _1})$ and in that case $n(\alpha $, $i) < 
(\beta $, $j)$ for all $n \in  {\open Z}$. Let $\hat{R} = K[[G]]$ be defined as in 
Definition 10, for this $G$ and for a field K of cardinality $\leq \hao$.
Define p-supp and truncation analogously to Definition 10.

 Let $Y_1$ be the smallest
valuation subring of $\hat{R} $ containing $X^{(\alpha,0)}$ for all
$\alpha  \in  \omega _1$ and let $r_\nu  = X^{(\nu , 0)}$ for $\nu  < \omega _1$.
Then $Q = \cup_{\nu < \goo} r_\nu^{-1}R$ and $\Gamma _{Y_1}(Q/Y_1) = 1$.

Define $G_i$ to be the subgroup of $G$ generated by $\{(\alpha $, $j)\colon \alpha  
< \omega _1$, $j < i\}$ and let $\hat{R} _i$ be the subring of $\hat{R} $ 
consisting of all $y \in  \hat{R} $ such that $\supp(y) \subseteq  G_i$. We 
shall now define a continuous chain $\{Y_i\colon i < 2^{\aleph _1}\}$ of valuation 
subrings of $R$ such that for each $i$, $Y_i \subseteq  \hat{R} _i$, and
is closed under truncation --- more precisely, if $y \in Y_i$ and $\gd
\in \lim(\goo)$, then $y\rest (\gd, 0) \in Y_i$.   Moreover, we  require that
$$\hbox{for all } i < j, Y_j \cap  \hat{R} _i = Y_i.\leqno(\hbox{\pounds})$$

We will let $R$ be the union of the $Y_i$. Then certainly $\Gamma _R(Q/R) = 1$,
since we will have added no limits of elements of $Y_1$. Moreover, we will do the
construction so that for each family $\{e^\tau _\sigma \colon \sigma  < \tau  < 
\omega _1\}$ as in (\dag ), there is a family $\{c_\sigma \colon \sigma  < 
\omega _1\}$ as in (\dag \dag ). There are only $(2^{\aleph _1})^{\aleph _1} = 
2^{\aleph _1}$ possible families $\{e^\tau _\sigma \colon \sigma  < \tau  <
\omega _1\}$ and the cofinality of $2^{\aleph _1} > \aleph _1$, so we can
arrange our construction such that each such family is considered at some stage
$i < 2^{\aleph _1}$ where  $\{e^\tau _\sigma \colon \sigma  < \tau  <
\omega _1\} \se Y_i.$

Suppose now that $Y_i$ has been constructed and that at this stage we are 
considering the family $\{e^\tau _\sigma \colon \sigma  < \tau  < \omega _1\} 
\subseteq  Y_i$. Our plan is to let $Y_{i+1}$ be the smallest valuation subring
of $\hat{R} $ which is closed under truncation and contains $Y_i \cup  \{X^{(\alpha , i)}\colon \alpha  <
\omega _1\}$  and a family $\{c_\sigma \colon \sigma  < \omega _1\}$ as in 
(\dag \dag ); we must choose $\{c_\sigma \colon \sigma  < \omega _1\}$ so that (\pounds) 
holds. We define the $c_\sigma $ by induction. Let $c_0 = 1$; if $c_\sigma $ 
has been defined, let $$c_{\sigma +1} = e^{\sigma +1}_\sigma c_\sigma  +
X^{(\sigma, i)}.$$ If $\delta $ is a limit ordinal, let $c_\delta $ be
the unique element of $\hat{R} $ with support in  $
\oplus _{j<i}\oplus _{\alpha <\omega _1}{\open Z}(\alpha $, $j)$ which 
represents the limit of $\langle e^\delta _\sigma c_\sigma \colon \sigma  < 
\delta \rangle $ in $\hat{R} /\cap _{\sigma <\delta }r_\sigma \hat{R}$.

We sketch the argument that this construction works. First, by induction
one can show that for all $\gt$, $(\gt, i) \in \psupp(c_\gs)$ if and
only if $\gt < \gs$. Secondly,
 the $c_\sigma $ are algebraically 
independent over the quotient field of $\hat{R} _i$.
More generally, if $\gs_1 < \ldots \gs_n < \gd \leq \gs_{n + 1}$, then
$\{c_{\gs_1}, \ldots, c_{\gs_n}, c_{\gs_{n + 1}}\rest\gd\}$ is
algebraically independent over $\hat{R}_i$. (This is proved
 by an
argument on supports.)   Now an arbitrary element of $Y_{i+1}$ has the form
$$\frac{f(c_{\gs_1}, \ldots, c_{\gs_n})\rest(\gd, 0)}{g(c_{\gs_1}, \ldots,
c_{\gs_n})\rest(\grg, 0)}$$ where $f$ and $g$ have coefficients in $Y_i$. By
an argument on supports, using the above facts, one shows that if this
element belongs to $Y_{i + 1} \cap \hat{R}_i$, then it belongs to $Y_i$.
                      \qed
\bigskip

\bigskip

\centerline{{\bf REFERENCES}}

\bigskip

\noindent
[BFS] S. Bazzoni, L. Fuchs, and L. Salce, {\it The hierarchy of
uniserial modules over a valuation domain}, preprint.

\bigskip

\noindent
[BS1] S. Bazzoni and L. Salce, {\it On non-standard uniserial modules over 
valuation domains and their quotients}, J. Algebra {\bf 128} (1990), 292-305.

\bigskip

\noindent
[BS2] S. Bazzoni and L. Salce, {\it Elongations of uniserial modules over 
valuation domains}, to appear in J. Algebra.

\bigskip

\noindent
[BS3] S. Bazzoni and L. Salce, {\it Equimorphism classes of uniserial modules 
over valuation domains}, to appear in Archiv. der Math.

\bigskip

\noindent
[D] K. Devlin, {\bf Constructibility}, Springer-Verlag (1984).

\bigskip

\noindent
[DSh] K. Devlin and S. Shelah, {\it A weak version of $\diamondsuit $ which follows from $2^{\aleph _0} <
2^{\aleph _1}$,} Israel J. Math. {\bf 29} (1978), 239--247.

\bigskip

\noindent
[E] P. Eklof, {\it A transfer theorem for non-standard uniserials},
 Proc. Amer. Math. Soc. {\bf 114} (1992), 593--600.

\bigskip

\noindent
[EM] P. Eklof and A. Mekler, {\bf Almost Free Modules}, North-Holland (1990).

\bigskip

\noindent
[FrG] B. Franzen and R. G\"obel, {\it Nonstandard uniserial modules over 
valuation domains,} Results in Math. {\bf 12} (1987), 86--94.

\bigskip

\noindent
[FS] L. Fuchs and L. Salce, {\bf Modules Over Valuation Domains}, Marcel 
Dekker (1985).

\bigskip

\noindent
[FSh] L. Fuchs and S. Shelah, {\it Kaplansky}'{\it s problem on valuation 
rings}, Proc. Amer. Math. Soc. {\bf 105} (1989), 25--30.

\bigskip

\noindent
[O1] B. L. Osofsky, {\it Constructing Nonstandard Uniserial Modules over
Valuation Domains},  {\bf Azumaya Algebras, Actions, and Modules},
Contemporary Mathematics {\bf 124} (1992),
 151--164.

\bigskip

\noindent
[O2] B. L. Osofsky, {\it A construction of nonstandard uniserial modules over
valuation domains},  Bull. Amer. Math. Society, {\bf 25} (1991), 89--97.

\bigskip

\noindent
[Sh] S. Shelah, {\it Nonstandard uniserial module over a uniserial domain 
exists}, Lecture Notes in Math vol. 1182, Springer-Verlag (1986), pp. 135--150.
\end{document}